\documentclass[aop,MSNbibl,seceqn,citesort,dvips]{arximspdf}

\doi{10.1214/11-AOP669}
\volume{40}
\issue{5}
\pubyear{2012}
\firstpage{2008}
\lastpage{2033}

\makeatletter
\newcommand{\EE}{\mathbb{E}}
\newcommand{\ZZ}{\mathbb{Z}}
\newcommand{\NN}{\mathbb{N}}
\newcommand{\1}{1}
\newcommand{\Exp}{\mathrm{E}}
\newcommand{\E}{\Exp}
\renewcommand{\Pr}{\mathrm{P}}
\newcommand{\xrightarrow}[1]{\stackrel{#1}{\rightarrow}}
\newcommand{\dto}{\xrightarrow{d}}
\newcommand{\wto}{\xrightarrow{w}}
\newcommand{\vto}{\xrightarrow{v}}
\newcommand{\fidi}{\xrightarrow{\mathrm{fidi}}}
\newcommand{\toi}{\to\infty}
\newcommand{\eind}{\stackrel{d}{=}}
\newcommand{\rmd}{\mathrm{d}}

\newcommand{\som}{\sum }
\renewcommand{\le}{\leq}
\renewcommand{\ge}{\geq}

 \newtheorem{theo}{Theorem}[section]
 \newtheorem{lem}[theo]{Lemma}
\newproclaim{ex}[theo]{Example}
\newproclaim{cond}[theo]{Condition}
\newcommand{\eqref}[1]{(\ref{#1})}
\newcommand{\fraca}[2]{{#1}/{#2}}
\makeatother

\begin{document}
\begin{frontmatter}

\title{A functional limit theorem for dependent sequences with
infinite variance stable limits}
\runtitle{Functional limit theorem}

\begin{aug}
\author[A]{\fnms{Bojan} \snm{Basrak}\corref{}\thanksref{a1}\ead[label=e1]{bbasrak@math.hr}},
\author[B]{\fnms{Danijel} \snm{Krizmani\'{c}}\ead[label=e2]{dkrizmanic@math.uniri.hr}}
\and
\author[C]{\fnms{Johan} \snm{Segers}\thanksref{a3}\ead[label=e3]{johan.segers@uclouvain.be}}
\runauthor{B. Basrak, D. Krizmani\'{c} and J. Segers}
\affiliation{University of Zagreb, University of Rijeka and
Universit\'e catholique~de~Louvain}
\address[A]{B. Basrak\\ Department of Mathematics\\
University of Zagreb\\ Bijeni\v cka 30, Zagreb\\ Croatia\\
\printead{e1}} 
\address[B]{D. Krizmani\'{c}\\
Department of Mathematics\\
University of Rijeka\\
Omladinska 14, Rijeka\\
Croatia\\\printead{e2}}
\address[C]{J. Segers\\
Institut de statistique, biostatistique et\\ \quad  sciences actuarielles\\
Universit\'e catholique de Louvain\\
Voie du Roman Pays 20\\
 B-1348 Louvain-la-Neuve\\
Belgium\\
\printead{e3}}
\end{aug}
\thankstext{a1}{Supported by the research Grant MZO\v S project no.
037-0372790-2800 of the Croatian government.}
\thankstext{a3}{Supported by IAP research network
Grant no. P6/03 of the Belgian government (Belgian Science Policy)
and from the contract ``Projet d'Actions de Recherche Concert\'ees''
no. 07/12/002 of the Communaut\'e fran\c{c}aise de Belgique,
granted by the Acad\'emie universitaire Louvain.}

\received{\smonth{1} \syear{2010}}
\revised{\smonth{11} \syear{2010}}

%
\begin{abstract}
Under an appropriate regular variation condition, the affinely
normalized partial sums of a sequence of independent and identically
distributed random variables converges weakly to a non-Gaussian
stable random variable. A functional version of this is known to be
true as well, the limit process being a stable L\'evy process. The
main result in the paper is that for a stationary, regularly varying
sequence for which clusters of high-threshold excesses can be broken
down into asymptotically independent blocks, the properly centered
partial sum process still converges to a stable L\'evy process. Due
to clustering, the L\'evy triple of the limit process can be
different from the one in the independent case. The convergence
takes place in the space of c\`adl\`ag functions endowed with
Skorohod's $M_1$ topology, the more usual $J_1$ topology being
inappropriate as the partial sum processes may exhibit rapid
successions of jumps within temporal clusters of large values,
collapsing in the limit to a single jump. The result rests on a new
limit theorem for point processes which is of independent interest.
The theory is applied to moving average processes, squared
$\operatorname{GARCH}(1,1)$ processes and stochastic volatility models.
\end{abstract}

%
\begin{keyword}[class=AMS]
\kwd[Primary ]{60F17}
\kwd{60G52}
\kwd[; secondary ]{60G55}
\kwd{60G70}.
\end{keyword}
\begin{keyword}
\kwd{Convergence in distribution}
\kwd{functional limit theorem}
\kwd{GARCH}
\kwd{mixing}
\kwd{moving average}
\kwd{partial sum}
\kwd{point processes}
\kwd{regular variation}
\kwd{stable processes}
\kwd{spectral processes}
\kwd{stochastic volatility}.
\end{keyword}

\end{frontmatter}

\section{Introduction}\label{intro}
\label{Sintro}

Consider a stationary sequence of random
variables $(X_n)_{n\geq 1}$ and its accompanying sequence of partial
sums $S_n=X_1+\cdots+X_n, {n\geq1}$. The main\vadjust{\goodbreak} goal of this paper
is to investigate the asymptotic distributional behavior of the
$D[0,1]$ valued process
\[
V_{n}(t) = a_{n}^{-1} \bigl(S_{\lfloor nt \rfloor} - \lfloor nt \rfloor
b_{n}\bigr), \qquad t
\in[0,1],
\]
under the properties of weak dependence and regular variation with
index $\alpha\in
(0,2)$, where
$(a_{n})_{n}$ is a sequence of positive real numbers such that
%
\begin{equation}\label{eniz}
n \Pr( |X_{1}| > a_{n}) \to1
\end{equation}
as $n \to\infty$  and
\[
b_{n} = \E \bigl( X_{1}   1_{\{ |X_{1}| \le a_{n} \}}  \bigr).
\]
Here, $\lfloor x \rfloor$ represents the integer part of the real
number $x$,
and $D[0, 1]$ is the space of real-valued c\`adl\`ag functions on
$[0, 1]$.

Recall that if the sequence $(X_n)_n$ is i.i.d. and if there exist
real sequences $a_n' > 0$ and $b_n'$ and a nondegenerate random
variable $S$ such that as $n\toi$
%
\begin{equation}
\label{ECLT}
\frac{S_n-b_n'}{a_n'} \dto S ,
\end{equation}
then $S$ is necessarily an $\alpha$-stable random variable. In
standard terminology, the law of $X_1$ belongs to the domain of
attraction of $S$. The domain of attraction of non-Gaussian stable
random variables can be completely characterized by an appropriate
regular variation condition; see \eqref{eregvar1} below. Classical
references in the i.i.d. case are the books by Gnedenko and
Kolmogorov~\cite{Gnedenko54}, Feller
\cite{Feller71} and Petrov~\cite{Petrov95}. In LePage et al.
\cite{LePage81} one can find an elegant probabilistic proof of
sufficiency and a nice representation of the limiting distribution.

Weakly dependent sequences can exhibit very similar behavior. The
first results in this direction were rooted in martingale theory (see
Durrett and Resnick~\cite{Durrett78}). In~\cite{Da83}, Davis
proved that if a regularly varying sequence of random
variables $(X_n)_n$ has tail index $0<\alpha<2$ and satisfies a strengthened
version of Leadbetter's $D$ and $D'$ conditions familiar from
extreme value theory, then \eqref{ECLT} holds for some
$\alpha$-stable random variable $S$ and properly chosen normalizing
sequences. These conditions are quite restrictive, however,
even excluding $m$-dependent sequences. Extensions to Davis's results
were provided in Denker and Jakubowski~\cite{DeJa89} and Jakubowski
and Kobus~\cite{JaKo89}, the latter paper being the first one in which
clustering of big values is allowed. Using classical blocking
techniques, necessary and sufficient conditions for convergence of sums
of weakly dependent random variables to stable laws are given in two
papers by Jakubowski~\cite{Ja93,Ja97}. The case of associated
sequences was treated in Dabrowski and Jakubowski~\cite{DaJa}.
In~\cite{DaHs95}, Davis and Hsing showed that sequences which satisfy a
regular variation condition for some $\alpha\in(0,2)$ and certain
mixing conditions also satisfy \eqref{ECLT} with an
$\alpha$-stable limit. Building upon the same point process
approach, Davis and Mikosch~\cite{DaMi98} generalized these results
to multivariate sequences. A survey of these results is to be found in
Bartkiewicz et
al.~\cite{BaJaMiWi09}, providing a detailed study of the conditions
for the convergence of the partial sums of a strictly stationary
process to an
infinite variance stable distribution. In this paper, the parameters of
the limiting distribution are determined in terms of some tail
characteristics of the underlying stationary sequence.

The asymptotic behavior of the processes $V_n$ as $n \to\infty$ is
an extensively studied subject in the probability literature, too. As
the index of regular variation $\alpha$ is assumed to be less than
$2$, the variance of $X_{1}$ is infinite. In the finite-variance
case, functional central limit theorems differ considerably and have
been investigated in greater depth (see, e.g.,
Billingsley~\cite{Billingsley68}, Herrndorf~\cite{Herrndorf85},
Merlev\`ede and Peligrad~\cite{Merlevede00}, and Peligrad and
Utev~\cite{Peligrad05}).

A functional limit theorem for the processes $V_n$ for infinite
variance i.i.d. regularly varying sequences $(X_n)$ was established in
Skorohod~\cite{Skorohod57}, a very readable proof of which can be
found in Resnick~\cite{Resnick07}. For stationary sequences, this
question was studied by Leadbetter and Rootz\'{e}n~\cite{Leadbetter88}
and Tyran-Kami{\'n}ska~\cite{TK10SPA}. Essentially, what they showed
is that the functional limit theorem holds in Skorohod's $J_1$ topology
if and only if certain point processes of extremes converge to a
Poisson random measure, which in turn is equivalent to a kind of
nonclustering property for extreme values. The implication is that for
many interesting models, convergence in the $J_1$ topology cannot hold.
This fact led Avram and Taqqu~\cite{Avram92} and Tyran-Kami{\'n}ska
\cite{TK10SPL} to opt for Skorohod's $M_1$ topology instead, a choice
which turns out to work for linear processes with regularly varying
innovations and nonnegative coefficients; see Section~\ref{Sflt} for
the definition of the $M_1$ topology.
For some more recent articles  with related but somewhat different
subjects we refer to Sly and Heyde~\cite{Sly08} who obtained
nonstandard limit theorems for functionals of regularly varying
sequences with long-range Gaussian dependence structure, and also to
Aue et al.~\cite{Aue08} who investigated the limit behavior of the
functional CUSUM statistic and its randomly permuted version for
i.i.d. random variables which are in the domain of attraction of a
strictly $\alpha$-stable law, for $\alpha\in(0,2)$.

The main theorem of our article shows that for a stationary,
regularly varying sequence for which clusters of high-threshold
excesses can be broken down into asymptotically independent blocks,
the properly centered partial sum process $(V_n(t))_{t\in[0,1]}$
converges to an $\alpha$-stable L\'evy process in the space
$D[0,1]$ endowed with Skorohod's $M_1$ metric under the condition
that all extremes within one such cluster have the same sign.
Our method of proof combines some ideas used in the
i.i.d. case by Resnick~\cite{Resnick86,Resnick07} with a new point
process convergence result and some particularities of the $M_1$
metric on $D[0,1]$ that can be found in Whitt~\cite{Whitt02}. The
theorem can be viewed as a generalization of results in Leadbetter
and Rootz\'{e}n~\cite{Leadbetter88} and Tyran-Kami{\'n}ska~\cite{TK10SPA}, where clustering of extremes is
essentially prohibited, and in Avram and Taqqu~\cite{Avram92} and
Tyran-Kami{\'n}ska~\cite{TK10SPL}, which are restricted to linear processes.

The paper is organized as follows. In Section~\ref{Sstatpoint} we
determine precise conditions needed to separate clusters of extremes
asymptotically. We also prove a new limit theorem for point
processes which is the basis for the rest of the paper and which is
of independent interest, too. In Section~\ref{Sflt} we state and
prove our main functional limit theorem. We also discuss possible
extensions of this result to other topologies. Finally, in
Section~\ref{Sexamples} several examples of stationary sequences
covered by our main theorem are discussed, in particular moving
averages and squared $\operatorname{GARCH}(1,1)$ processes.

%
\section{Stationary regularly varying sequences}
\label{Sstatpoint}

The extremal dynamics of a~regularly varying stationary time series
can be captured by its tail process, which is the conditional
distribution of the series, given that at a certain moment, it is far
away from the origin (Section~\ref{SSstatpointtail}). In
particular, the tail process allows explicit descriptions of the
limit distributions of various point processes of extremes
(Section~\ref{SSstatpointpoint}). The main result in this
section is Theorem~\ref{Tpointprocesscomplete}, providing the weak
limit of a sequence of time-space point processes, recording both
the occurrence times and the values of extreme values.

\subsection{Tail processes}
\label{SSstatpointtail}

Denote $\EE=\overline{\mathbb{R}} \setminus\{ 0 \}$ where
$\overline{\mathbb{R}}=[-\infty, \infty]$. The space $\EE$ is
equipped with the topology which makes it homeomorphic to $[-1, 1]
\setminus\{0\}$ (Euclidean topology) in the obvious way. In
particular, a set $B \subset\EE$ has compact closure if and only if
it is bounded away from zero, that is, if there exists $u > 0$ such
that $B \subset\EE_u = \EE\setminus[-u, u]$. Denote by
$C_{K}^{+}(\EE)$ the class of all nonnegative, continuous functions
on $\EE$ with compact support.

We say that a strictly stationary process $(X_{n})_{n \in
\mathbb{Z}}$ is \textit{(jointly) regularly varying} with index
$\alpha\in(0,\infty)$ if for any nonnegative integer $k$, the
$k$-dimensional random vector $\mathbf{X} = (X_{1}, \ldots,
X_{k})$ is multivariate regularly varying with index $\alpha$; that is,
for some (and then for every) norm $\|   \cdot  \|$ on
$\mathbb{R}^{k}$ there exists a random vector $\bolds{\Theta}$
on the unit sphere $\mathbb{S}^{k-1} = \{ x \in\mathbb{R}^{k} \dvtx
\|x\|=1 \}$ such that for every $u \in(0,\infty)$ and as $x \to
\infty$,
%
\begin{equation}\label{eregvar1}
\frac{\Pr(\|\mathbf{X}\| > ux, \mathbf{X} / \| \mathbf
{X} \| \in\cdot  )}{\Pr(\| \mathbf{X} \| >x)}
\wto u^{-\alpha} \Pr( \bolds{\Theta} \in\cdot ),
\end{equation}
the arrow ``$\wto$'' denoting weak convergence of finite measures.
For an extensive and highly-readable account of (multivariate)
regular variation, see the monograph by Resnick~\cite{Resnick07}.

Theorem 2.1 in Basrak and Segers~\cite{BaSe} provides a convenient
characterization of joint regular variation: it is necessary and
sufficient that there exists a process $(Y_n)_{n \in\mathbb{Z}}$
with $\Pr(|Y_0| > y) = y^{-\alpha}$ for $y \ge1$ such that as $x
\to\infty$,
%
\begin{equation}\label{etailprocess}
 \bigl( (x^{-1}X_n)_{n \in\ZZ}    \mid   |X_0| > x  \bigr)
\fidi(Y_n)_{n \in\ZZ},
\end{equation}
where ``$\fidi$'' denotes convergence of finite-dimensional
distributions. The process $(Y_{n})_{n \in\mathbb{Z}}$ is called
the \textit{tail process} of $(X_{n})_{n \in\mathbb{Z}}$. Writing
$\Theta_n = Y_n / |Y_0|$ for $n \in\ZZ$, we also have
\[
 \bigl( (|X_0|^{-1}X_n)_{n \in\ZZ}    \mid   |X_0| > x  \bigr)
\fidi(\Theta_n)_{n \in\ZZ}
\]
(see Corollary 3.2 in~\cite{BaSe}). The process $(\Theta_n)_{n \in
\ZZ}$ is independent of $|Y_0|$ and is called the \textit{spectral
(tail) process} of $(X_n)_{n \in\ZZ}$. The law of $\Theta_0 = Y_0 /
|Y_0| \in\mathbb{S}^{0} = \{-1, 1\}$ is the spectral measure of the
common marginal distribution of the random variables~$X_i$. Regular
variation of this marginal distribution can be expressed in terms of
vague convergence of measures on $\EE$: for $a_n$ as in
\eqref{eniz} and as $n \to\infty$,
%
\begin{equation}
\label{eonedimregvar}
n \Pr( a_n^{-1} X_i \in\cdot  ) \vto\mu(   \cdot ),
\end{equation}
the Radon measure $\mu$ on $\EE$ being given by
%
\begin{equation}
\label{Emu}
\mu(\rmd x) = \bigl ( p   1_{(0, \infty)}(x) + q   1_{(-\infty,
0)}(x)  \bigr)   \alpha|x|^{-\alpha-1} \, \rmd x,
\end{equation}
where
\begin{eqnarray*}
p &=& \Pr(\Theta_0 = +1) = \lim_{x \to\infty} \frac{\Pr(X_i >
x)}{\Pr(|X_i| > x)}, \\
q &=& \Pr(\Theta_0 = -1) = \lim_{x \to\infty} \frac{\Pr(X_i <
-x)}{\Pr(|X_i| > x)}.
\end{eqnarray*}

\subsection{Point process convergence}
\label{SSstatpointpoint}

Define the time-space point processes
%
\begin{equation}
\label{Eppspacetime}
N_{n} = \sum_{i=1}^{n} \delta_{(i / n, X_{i} / a_{n})} \qquad\mbox
{for all $n\in\NN$}
\end{equation}
with $a_n$ as in \eqref{eniz}. The aim of this section is to
establish weak convergence of $N_n$ in the state space $[0, 1]
\times\EE_u$ for $u > 0$, where {$\EE_u = [-\infty, -u) \cup(u,
\infty]$}.
The limit process is a Poisson superposition of cluster processes,
whose distribution is determined by the tail process $(Y_i)_{i \in
\ZZ}$. Convergence of $N_n$ was already alluded to without proof in
Davis and Hsing~\cite{DaHs95} with a reference to
Mori~\cite{Mori77}.

To control the dependence in the sequence $(X_n)_{n \in\mathbb{Z}}$
we first have to assume that clusters of large values of $|X_{n}|$
do not last for too long.

\begin{cond}[(Finite mean cluster size)]
\label{cfinite-mean-cluster-size} There exists a positive integer sequence
$(r_{n})_{n \in\mathbb{N}}$ such that $r_{n} \to\infty$ and
$r_{n} / n \to0$ as $n \to\infty$ and such that for every $u > 0$,
%
\begin{equation}
\label{eanticluster}
\lim_{m \to\infty} \limsup_{n \to\infty}
\Pr \Bigl( \max_{m \le|i| \le r_{n}} |X_{i}| > ua_{n} \bigm|
|X_{0}|>ua_{n}  \Bigr) = 0.
\end{equation}
\end{cond}

Put $M_{1,n} = \max\{ |X_{i}| \dvtx  i=1, \ldots, n \}$ for $n \in
\NN$. In Proposition 4.2 in~\cite{BaSe}, it has been shown that
under  Condition~\ref{cfinite-mean-cluster-size} we have $\theta>
0$, where
%
\begin{eqnarray}
\label{Ethetaspectral}
\theta&=& \lim_{r \to\infty} \lim_{x \to\infty} \Pr
(M_{1,r} \le x    \mid   |X_{0}|>x  ) \nonumber
\\[-8pt]
\\[-8pt]
&=& \Pr\Bigl({\sup_{i\ge1}} |Y_{i}| \le1\Bigr) = \Pr\Bigl({
\sup_{i\le-1}} |Y_{i}| \le1\Bigr).
\nonumber
\end{eqnarray}
%
Moreover $\Pr( \lim_{|n| \to\infty} | Y_n | =
0 ) = 1$, and, for every $u \in(0, \infty)$ and as $n \to\infty$,
%
\begin{equation}
\label{Erunsblocks}
\Pr(M_{1,r_n} \leq a_n u \mid| X_0 | > a_n u)
= \frac{\Pr(M_{1,r_n} > a_n u)}{r_n \Pr(| X_0 | > a_n u)} + o(1)
\to\theta,
\end{equation}
and thus
\[
\lim_{n \to\infty} \E\Biggl [ \sum_{t=1}^{r_n} \1_{(a_n u, \infty
)}(X_t)    \Bigm|   M_{1,r_n} > a_n u  \Biggr] = \frac{1}{\theta} <
\infty.
\]
This explains why we call Condition~\ref{cfinite-mean-cluster-size}
the ``finite mean cluster size'' condition. In the setting of
Theorem~\ref{Tpointprocesscomplete} below, the
quantity $\theta$ in \eqref{Ethetaspectral} is the \textit{extremal
index} of the sequence $(| X_n |)_{n \in\mathbb{Z}}$ (see~\cite{BaSe}, Remark 4.7): for all $u \in(0, \infty)$ and as $n \to\infty$,
\[
\Pr(M_{1,n} \le a_n u)
= \bigl ( \Pr(| X_1 | \le a_n u)  \bigr)^{n \theta} + o(1)
\to e^{-\theta u^{-\alpha}}.
\]
See Section~\ref{SSSfmcs} for further discussion.

Since $\Pr(M_{1,r_n} > a_n u) \to0$ as $n \to\infty$, we call the
point process
\[
\sum_{i=1}^{r_n} \delta_{(a_n u)^{-1} X_i} \qquad\mbox{conditionally
on }   M_{1,r_n} > a_n u
\]
a \textit{cluster process}, to be thought of as a cluster of
exceptionally large values occurring in a relatively short time
span. Theorem 4.3 in~\cite{BaSe} yields the weak convergence of the
sequence of cluster processes in the state space $\EE$,
%
\begin{equation}
\label{Eclusterprocess}
 \Biggl( \sum_{i=1}^{r_n} \delta_{(a_n u)^{-1} X_i}    \Bigm|
M_{1,r_n} > a_n u  \Biggr)
\dto \biggl( \sum_{n \in\mathbb{Z}} \delta_{Y_n}    \Bigm|
\sup_{i \le-1} |Y_i| \le1  \biggr).
\end{equation}
Note that since $|Y_n| \to0$ almost surely as $|n| \to\infty$, the
point process $\sum_n \delta_{Y_n}$ is well defined in $\EE$. By
\eqref{Ethetaspectral}, the probability of the conditioning event
on the right-hand side of \eqref{Eclusterprocess} is nonzero.

To establish convergence of $N_n$ in \eqref{Eppspacetime}, we need
to impose a certain mixing condition denoted by $\mathcal{A}'(a_n)$
which is slightly stronger than the condition $\mathcal{A}(a_n)$
introduced in Davis and Hsing~\cite{DaHs95}.

\begin{cond}[($\mathcal{A}'(a_{n})$)]
\label{cmixcond} There exists a sequence of positive integers
$(r_{n})_{n}$ such that $r_{n} \to\infty$ and $r_{n} / n \to0$ as
$n \to\infty$ and such that for every\vadjust{\goodbreak} $f \in C_{K}^{+}([0,1] \times
\EE)$, denoting $k_{n} = \lfloor n / r_{n} \rfloor$, as $n \to
\infty$,
%
\begin{equation}\label{emixcon}
  \Exp \Biggl[ \exp \Biggl\{ - \sum_{i=1}^{n} f \biggl (\frac{i}{n},
\frac{X_{i}}{a_{n}}
 \biggr)  \Biggr\}  \Biggr]
- \prod_{k=1}^{k_{n}} \Exp \Biggl[ \exp \Biggl\{ - \sum
_{i=1}^{r_{n}} f  \biggl(\frac{kr_{n}}{n}, \frac{X_{i}}{a_{n}}
\biggr)  \Biggr\}  \Biggr] \to0.\hspace*{-40pt}
\end{equation}
\end{cond}

It can be shown that Condition~\ref{cmixcond} is implied by the
strong mixing property (see Krizmani\'c~\cite{Kr10}). {Recall $\EE_u
= \EE\setminus[-u,u]$.}

\begin{theo}
\label{Tpointprocesscomplete} If Conditions~\ref{cfinite-mean-cluster-size}
and~\ref{cmixcond} hold, then for every $u \in(0, \infty)$ and as
$n \to\infty$,
\[
N_n { |_{[0, 1] \times\EE_u} } \dto N^{(u)}
= \sum_i \sum_j \delta_{(T^{(u)}_i, u Z_{ij})}  \big|_{[0, 1]
\times\EE_u}
\]
in $[0, 1] \times\EE_u$, where:
\begin{longlist}[(2)]
\item[(1)]$\sum_i \delta_{T^{(u)}_i}$ is a homogeneous Poisson process on
$[0, 1]$ with intensity~$\theta u^{-\alpha}$;
\item[(2)]$(\sum_j \delta_{Z_{ij}})_i$ is an i.i.d. sequence of point
processes in $\EE$, independent of $\sum_i \delta_{T^{(u)}_i}$, and
with common distribution equal to the weak limit in \eqref{Eclusterprocess}.
\end{longlist}
\end{theo}

It can be shown that Theorem~\ref{Tpointprocesscomplete} is still
valid if $\EE_u$ is
replaced by {$\overline{\EE_u} = [-\infty, -u] \cup[u, \infty]$}.

\begin{pf*}{Proof of Theorem~\ref{Tpointprocesscomplete}}
Let $(X_{k,j})_{j \in\NN}$, with $k \in\NN$, be independent copies
of $(X_j)_{j \in\NN}$, and define
\[
\hat{N}_n
= \sum_{k=1}^{k_n} \hat{N}_{n,k}
\qquad\mbox{with }
\hat{N}_{n,k} = \sum_{j=1}^{r_n}
\delta_{(k r_n/n, X_{k,j}/a_n)}.
\]
By Condition~\ref{cmixcond}, the weak limits of $N_n$ and
$\hat{N}_n$ must coincide. By Kallenberg~\cite{Kallenberg83},
Theorem 4.2, it is enough to show that the Laplace functionals
of~$\hat{N}_n$ converge to those of $N^{(u)}$. Take $f \in
C_K^+([0,1] \times\EE_u)$. We extend $f$ to the whole of $[0, 1]
\times\EE$ by setting $f(t,x)=0$ whenever $|x| \le u$; in this way,
$f$ becomes a bounded, nonnegative and continuous function on $[0,
1] \times\EE$. There exists $M \in(0, \infty)$ such that $0 \le
f(t, x) \le M   1_{[-u,u]^c}(x)$. Hence as $n \to\infty$,
\begin{eqnarray*}
1
&\ge&\Exp e^{-\hat{N}_{n,k} f}
 \ge \Exp e^{-M\sum_{i=1}^{r_n} 1(| X_i | > a_n u)} \\
&\ge&1 - M r_n \Pr(| X_0 | > a_n u) = 1 - O(k_n^{-1}).
\end{eqnarray*}
In combination with the elementary bound $0 \le- \log z - (1-z) \le
(1-z)^2/z$ for $z \in(0,1]$, it follows that as $n \to\infty$,
\[
-\log\E e^{- \hat{N}_{n} f}
= - \sum_{k=1}^{k_n} \log\E e^{- \hat{N}_{n,k} f}
= \sum_{k=1}^{k_n} (1 - \E e^{- \hat{N}_{n,k} f}) + O(k_n^{-1}).
\]

By \eqref{Erunsblocks}, $k_n \Pr(M_{1,r_n} >a_n u ) \to\theta
u^{-\alpha}$ for $u \in(0, \infty)$ and as $n \to\infty$. Hence
%
\begin{eqnarray}\label{EKall1}
&& \quad \sum_{k=1}^{k_n} (1 - \E e^{- \hat{N}_{n,k} f}) \nonumber\\
&& \quad  \qquad =
k_n \Pr(M_{1,r_n}>a_n u )
\frac{1}{k_n} \sum_{k=1}^{k_n}
\Exp \bigl[ 1- e^{- \sum_{j=1}^{r_n} f(k r_n/n, X_{j}/a_n) }
|
M_{1,r_n}>a_n u  \bigr]\\
&& \quad  \qquad =
\theta u^{-\alpha}   \frac{1}{k_n} \sum_{k=1}^{k_n}
\Exp\bigl [ 1- e^{- \sum_{j=1}^{r_n} f(k r_n/n, X_j/a_n)}
  |  M_{1,r_n}>a_n u  \bigr] + o(1).\nonumber
\end{eqnarray}
Let the random variable $T_n$ be uniformly distributed on $\{k r_n /
n \dvtx  k = 1, \ldots, k_n\}$ and independent of $(X_j)_{j \in
\mathbb{Z}}$. By the previous display, as $n \to\infty$,
\[
\sum_{k=1}^{k_n} (1 - \Exp e^{- \hat{N}_{n,k} f})
= \theta u^{-\alpha} \Exp \bigl[ 1- e^{- \sum_{j=1}^{r_n}
f(T_n, u X_j/(ua_n)) }   |  M_{1,r_n}>a_n u  \bigr] + o(1).
\]
The sequence $T_n$ converges in law to a uniformly distributed
random variable~$T$ on $(0,1)$. By \eqref{Eclusterprocess} and by
independence of sequences $(T_n)$ and $(X_n)$
\[
 \Biggl( T_n, \sum_{i=1}^{r_n} \delta_{a_n^{-1} X_i}    \Bigm|
M_{1,r_n} > a_n u  \Biggr)
\dto \biggl(T, \sum_{n \in\mathbb{Z}} \delta_{u Z_n}  \biggr),
\]
where $\sum_n \delta_{Z_{n}}$ is a point process on $\EE$,
independent of the random variable $T$ and with distribution equal
to the weak limit in \eqref{Eclusterprocess}.
Thus the expressions in \eqref{EKall1}
converge as $n \to\infty$ to
%
\begin{equation}
\label{Ecompletelimit}
\theta u^{-\alpha}
\Exp \bigl[ 1- e^{- \sum_j f(T, u Z_j) }  \bigr]
= \int_0^1 \Exp \bigl[ 1 - e^{-\sum_j f(t, u Z_j) }  \bigr]
\theta u^{-\alpha} \, \rmd t.
\end{equation}
It remains to be shown that the right-hand side above equals $-\log
\Exp e^{-N^{(u)}f}$ for $N^{(u)}$ as in the theorem.

Define $g(t) = \Exp\exp\{ - \sum_j f(t, u Z_j) \}$ for $t \in[0,
1]$. Since $\sum_i \delta_{T_i^{(u)}}$ is independent of the i.i.d.
sequence $(\sum_j \delta_{Z_{ij}})_i$,
\begin{eqnarray*}
\Exp e^{- N^{(u)}f}
&=& \Exp e^{- \sum_i \sum_j f(T^{(u)}_i, u Z_{ij})} \\
&=& \Exp \biggl[ \prod_{i} \Exp \bigl( e^{- \sum_j f(T^{(u)}_i, u
Z_{ij})}   \mid\bigl (T^{(u)}_k\bigr)_k  \bigr)  \biggr]
= \Exp e^{\sum_i \log g(T^{(u)}_i)}.
\end{eqnarray*}
The right-hand side is the Laplace functional of a homogeneous
Poisson process on $[0,1]$ with intensity $\theta u^{-\alpha}$
evaluated in the function $- \log g$. Therefore, it is equal to
\[
\exp \biggl( - \int_0^1 \{1 - g(t)\} \theta u^{-\alpha} \, \rmd t
 \biggr)
\]
(see, e.g., Embrechts et al.~\cite{Em97}, Lemma 5.1.12; note
that $0 \le g \le1$). By the definition of $g$, the integral in the
exponent is equal to the one in \eqref{Ecompletelimit}. This
completes the proof of the theorem.
\end{pf*}

%
\section{Functional limit theorem}
\label{Sflt}

The main result in the paper states convergence of the partial sum
process $V_n$ to a stable L\'evy process in the space $D[0, 1]$
equipped with Skorohod's $M_1$ topology. The core of the proof rests
on an application of the continuous mapping theorem: the partial sum
process~$V_n$ is represented as the image of the time-space point
process $N_n$ in~\eqref{Eppspacetime} under a certain summation
functional. This summation functional enjoys the right continuity
properties by which the weak convergence of $N_n$ in
Theorem~\ref{Tpointprocesscomplete} transfers to weak convergence
of $V_n$.

The definition and basic properties of the $M_1$ topology are
recalled in Section~\ref{SSM1}. In
Section~\ref{SSsumfunct}, the focus is on the summation
functional and its continuity properties. The main result of the
paper then comes in Section~\ref{SSmain}. The conditions
entering this theorem are discussed in Section~\ref{SSdisc},
while Section~\ref{SSsimpl} provides some simplifications.


\subsection{The $M_1$ topology}
\label{SSM1}

The metric $d_{M_{1}}$ that generates the $M_{1}$ topology on $D[0,
1]$ is defined using completed graphs. For $x \in D[0,1]$ the
\textit{completed graph} of $x$ is the set
\[
\Gamma_{x}
= \{ (t,z) \in[0,1] \times\mathbb{R} \dvtx  z= \lambda x(t-) + (1-\lambda
)x(t)  \mbox{ for some } \lambda\in[0,1] \},
\]
where $x(t-)$ is the left limit of $x$ at $t$. Besides the points of
the graph $ \{ (t,x(t)) \dvtx  t \in[0,1] \}$, the completed graph of
$x$ also contains the vertical line segments joining $(t,x(t))$ and
$(t,x(t-))$ for all discontinuity points $t$ of $x$. We define an
\textit{order} on the graph $\Gamma_{x}$ by saying that $(t_{1},z_{1})
\le(t_{2},z_{2})$ if either (i) $t_{1} < t_{2}$ or (ii) $t_{1} =
t_{2}$ and $|x(t_{1}-) - z_{1}| \le|x(t_{2}-) - z_{2}|$. A
\textit{parametric representation} of the completed graph $\Gamma_{x}$
is a continuous nondecreasing function $(r,u)$ mapping $[0,1]$ onto
$\Gamma_{x}$, with $r$ being the time component and $u$ being the
spatial component. Let $\Pi(x)$ denote the set of parametric
representations of the graph $\Gamma_{x}$. For $x_{1},x_{2} \in
D[0,1]$ define
\[
d_{M_{1}}(x_{1},x_{2})
= \inf\bigl\{ \|r_{1}-r_{2}\|_{[0,1]} \vee\|u_{1}-u_{2}\|_{[0,1]} \dvtx
(r_{i},u_{i}) \in\Pi(x_{i}), i=1,2 \bigr\},
\]
where $\|x\|_{[0,1]} = \sup\{ |x(t)| \dvtx  t \in[0,1] \}$. This
definition introduces $d_{M_{1}}$ as a~metric on $D[0,1]$. The
induced topology is called Skorohod's $M_{1}$ topology and is weaker
than the more frequently used $J_{1}$ topology which is also due to
Skorohod.

The $M_{1}$ topology allows for a jump in the limit function $x \in
D[0, 1]$ to be approached by multiple jumps in the converging
functions $x_n \in D[0, 1]$. Let, for instance,
\[
x_n(t)  = \tfrac{1}{2} 1_{[\fraca{1}{2}-\fraca{1}{n}, \fraca
{1}{2})}(t) + 1_{[\fraca{1}{2}, 1]}(t),  \qquad
x(t)  = 1_{[\fraca{1}{2}, 1]}(t)
\]
for $n \geq3$ and $t \in[0, 1]$. Then $d_{M_{1}}(x_{n},x) \to0$
as $n \to\infty$, although $(x_{n})_n$ does not converge to $x$ in
either the uniform or the $J_{1}$ metric. For more discussion of the
$M_{1}$ topology we refer to Avram and Taqqu~\cite{Avram92} and
Whitt~\cite{Whitt02}.

\subsection{Summation functional}
\label{SSsumfunct}

Fix $0 < v < u < \infty$. The proof of our main theorem depends on
the continuity properties of the summation functional
\[
\psi^{(u)} \dvtx \mathbf{M}_{p}([0,1] \times\EE_{v}) \to D[0,1]
\]
defined by
\[
\psi^{(u)}  \biggl( \som_{i}\delta_{(t_{i}, x_{i})}  \biggr) (t)
= \sum_{t_{i} \le t} x_{i}  1_{\{u < |x_i| < \infty\}}, \qquad t \in
[0, 1].
\]
Observe that $\psi^{(u)}$ is well defined because $[0,1] \times
\EE_{u}$ is a relatively compact subset of $[0,1] \times\EE_{v}$.
The space $\mathbf{M}_p$ of Radon point measures is equipped with
the vague topology, and $D[0, 1]$ is equipped with the $M_1$
topology.

We will show that $\psi^{(u)}$ is continuous on the set $\Lambda=
\Lambda_{1} \cap\Lambda_{2}$, where
\begin{eqnarray*}
\Lambda_{1} &=&
\bigl\{ \eta\in\mathbf{M}_{p}([0,1] \times\EE_{v}) \dvtx
\eta( \{0,1 \} \times\EE_{u}) = 0 = \eta([0,1] \times\{ \pm\infty
, \pm u \}) \bigr\}, \\
 \Lambda_{2}& =&
\bigl\{ \eta\in\mathbf{M}_{p}([0,1] \times\EE_{v}) \dvtx
\eta\bigl( \{ t \} \times(v, \infty]\bigr) \wedge\eta\bigl( \{ t \} \times
[-\infty,-v)\bigr) = 0
\\
&&\hspace*{220pt}\mbox{for all $t \in[0,1]$} \bigr\};
\end{eqnarray*}
we write $s \wedge t$ for $\min(s, t)$. Observe that the elements
of $\Lambda_2$ have the property that atoms with the same time
coordinate are all on the same side of the time axis.

\begin{lem}
\label{lprob1} Assume that with probability one, the tail process
$(Y_{i})_{i \in\mathbb{Z}}$ in \eqref{etailprocess} has no two
values of the opposite sign. Then $ \Pr( N^{(v)} \in\Lambda) =
1$.
\end{lem}

The assumption that the tail process cannot switch sign will reappear
in our main result, Theorem~\ref{t2}. For linear processes, for
instance, it holds as soon as all coefficients are of the same sign.

\begin{pf}
From the definition of the tail process $(Y_{i})_{i \in\mathbb{Z}}$
we know that $\Pr(Y_{i}= \pm\infty)=0$ for any $i \in\mathbb{Z}$.
Moreover, by the spectral decomposition $Y_i = |Y_0| \Theta_i$ into
independent components $|Y_0|$ and $\Theta_i$ with $|Y_0|$ a Pareto
random variable, it follows that $Y_i$ cannot have any atoms
except possibly at the origin. As a consequence, it holds with probability one
that $\sum_j \delta_{ v Y_{j}} (\{\pm u\} ) =0$ and thus that
$\sum_j \delta_{ v Z_{ij}} (\{\pm u\} ) =0$ as well. Together with
the fact that $\Pr( \sum_{i}\delta_{T_{i}^{(v)}} (\{0,1\}) = 0 )=1$
this implies\vspace*{-1pt} $\Pr( N^{(v)} \in\Lambda_{1})=1$.

Second, the assumption that with probability one the tail process
$(Y_{i})_{i \in\mathbb{Z}}$ has no two values of the opposite sign
yields $\Pr(N^{(v)} \in\Lambda_{2})=1$.
\end{pf}

\begin{lem}
\label{Lcontsf} The summation functional $\psi^{(u)} \dvtx
\mathbf{M}_{p}([0,1] \times\EE_{v}) \to D[0,1]$ is continuous on
the set $\Lambda$, when $D[0,1]$ is endowed with Skorohod's $M_{1}$
metric.\vadjust{\goodbreak}
\end{lem}

\begin{pf}
Suppose that $\eta_{n} \vto\eta$ in $\mathbf{M}_p$ for some $\eta
\in\Lambda$. We will show that $\psi^{(u)}(\eta_n) \to
\psi^{(u)}(\eta)$ in $D[0, 1]$ according to the $M_1$ topology. By
Whitt ~\cite{Whitt02}, Corollary 12.5.1, $M_1$ convergence for
monotone functions amounts to pointwise convergence in a dense
subset of points plus convergence at the endpoints. Our proof is
based on an extension of this criterion to piecewise monotone
functions. This cut-and-paste approach is justified in view of
\cite{Whitt02}, Lemma 12.9.2, provided that the limit
function is continuous at the cutting points.

As $[0, 1] \times\EE_u$ is relatively compact in $[0, 1] \times
\EE_v$ there exists a nonnegative integer $k=k(\eta)$ such that
\[
\eta( [0,1] \times\EE_{u}) = k < \infty.
\]
By assumption, $\eta$ does not have any atoms on the horizontal
lines at $u$ or~$-u$. As a consequence, by Resnick ~\cite{Resnick07},
Lemma 7.1, there exists a positive integer~$n_{0}$
such that for all $n \ge n_{0}$ it holds that
\[
\eta_{n} ( [0,1] \times\EE_{u} ) = k.
\]
If $k = 0$, there is nothing to prove, so assume $k \ge1$, and let
$(t_i, x_i)$ for $i \in\{1, \ldots, k\}$ be the atoms of $\eta$ in
$[0,1] \times{\EE}_u$. By the same lemma, the~$k$ atoms
$(t_{i}^{(n)}, x_{i}^{(n)})$ of $\eta_n$ in $[0, 1] \times\EE_u$
(for $n \ge n_0$) can be labeled in such a~way that for $i \in\{1,
\ldots, k\}$, we have
\[
\bigl(t_{i}^{(n)},x_{i}^{(n)}\bigr) \to(t_{i},x_{i})  \qquad\mbox{as }  n \to
\infty.
\]
In particular, for any $\delta>0$ we can find a positive integer
$n_{\delta}$ such that for all $n \ge n_{\delta}$,
%
\begin{eqnarray}
\label{epointsdif}
\eta_{n} ([0,1] \times\EE_{u}) &=& k,\nonumber
\\[-8pt]
\\[-8pt]
\bigl|t_{i}^{(n)} - t_{i}\bigr| &<& \delta\quad\mbox{and} \quad
\bigl|x_{i}^{(n)} - x_{i}\bigr| < \delta  \qquad\mbox{for }  i=1, \ldots, k.
\nonumber
\end{eqnarray}
Let the sequence
\[
0< \tau_{1} < \tau_{2} <\cdots < \tau_{p} < 1
\]
be such that the sets $\{ \tau_{1}, \ldots, \tau_{p} \}$ and
$\{t_{1}, \ldots, t_{k} \}$ coincide. Note that \mbox{$p \le k$} always
holds, but since $\eta$ can have several atoms with the same time
coordinate, equality does not hold in general. Put $\tau_{0}=0$,
$\tau_{p+1}=1$, and take
\[
0 < r < \frac{1}{2}\min_{0 \le i \le p}|\tau_{i+1} -
\tau_{i}|.
\]
For any $t \in[0,1] \setminus\{ \tau_{1}, \ldots,
\tau_{p} \}$ we can find $\delta\in(0,u)$ such that
\[
\delta< r \quad\mbox{and} \quad\delta< \min_{1 \le i \le
p} |t - \tau_{i}|.
\]
Then relation \eqref{epointsdif}, for $n \ge n_{\delta}$,
implies that $t_{i}^{(n)} \le t$ is equivalent to $t_{i} \le
t$, and we obtain
\[
\bigl|\psi^{(u)}(\eta_{n})(t) - \psi^{(u)}(\eta)(t)\bigr| =  \biggl| \sum
_{t_{i}^{(n)} \le
t}x_{i}^{(n)} - \sum_{t_{i} \le
t}x_{i}  \biggr| \le\sum_{t_{i} \le t}\delta\le k\delta.
\]
Therefore
\[
\lim_{n \to\infty} \bigl|\psi^{(u)}(\eta_{n})(t) - \psi^{(u)}(\eta
)(t)\bigr| \le k\delta,
\]
and if we let $\delta\to0$, it follows that
$\psi^{(u)}(\eta_{n})(t) \to\psi^{(u)}(\eta)(t)$ as $n \to
\infty$.
Put
\[
v_i = \tau_i + r, \qquad i \in\{1, \ldots, p\}.
\]
For any $\delta< u \wedge r$, relation \eqref{epointsdif} and the
fact that $\eta\in\Lambda$ imply that the functions
$\psi^{(u)}(\eta)$ and $\psi^{(u)}(\eta_{n})$ ($n \ge n_{\delta}$)
are monotone on each of the intervals $[0,v_{1}], [v_{1},v_{2}],
\ldots, [v_{p},1]$. A combination of Corollary 12.5.1 and Lemma~12.9.2 in~\cite{Whitt02} yields
$d_{M_{1}}(\psi^{(u)}(\eta_{n}), \psi^{(u)}(\eta)) \to0$ as $n \to
\infty$. The application of Lem\-ma~12.9.2 is justified by continuity
of $\psi^{(u)}(\eta)$ in the boundary points $v_{1}, \ldots,
v_{p}$. We conclude that $\psi^{(u)}$ is continuous at $\eta$.
\end{pf}

\subsection{Main theorem}
\label{SSmain}

Let $(X_n)_n$ be a strictly stationary sequence of random variables,
jointly regularly varying with index $\alpha\in(0, 2)$ and tail
process $(Y_i)_{i \in\ZZ}$. The theorem below gives conditions
under which its partial sum process satisfies a nonstandard
functional limit theorem with a non-Gaussian $\alpha$-stable
L\'{e}vy process as a limit. Recall that the distribution of a
L\'{e}vy process $V( \cdot )$ is characterized by its
\textit{characteristic triple}, that is, the characteristic triple of the
infinitely divisible distribution of $V(1)$. The characteristic
function of $V(1)$ and the characteristic triple $(a, \nu, b)$ are
related in the following way:
\[
\E\bigl[e^{izV(1)}\bigr] = \exp \biggl( -\frac{1}{2}az^{2} + ibz + \int
_{\mathbb{R}}  \bigl( e^{izx}-1-izx 1_{[-1,1]}(x)  \bigr) \nu(\rmd
x)  \biggr)
\]
for $z \in\mathbb{R}$; here $a \ge0$, $b \in\mathbb{R}$ are
constants, and $\nu$ is a measure on $\mathbb{R}$ satisfying
\[
\nu( \{0\})=0 \quad\mbox{and} \quad\int_{\mathbb{R}}(|x|^{2}
\wedge1) \nu(\rmd x) < \infty;
\]
that is, $\nu$ is a L\'{e}vy measure. For a textbook treatment of
L\'{e}vy processes we refer to Bertoin~\cite{Bertoin96} and Sato~\cite{Sato99}. The description of the L\'{e}vy triple of the
limit process will be in terms of the measures $\nu^{(u)}$ ($u > 0$)
on $\EE$ defined for $x > 0$ by
%
\begin{eqnarray}
\label{Enuu}
\nu^{(u)}(x, \infty) &=&  u^{-\alpha}   \Pr \biggl( u
\sum_{i \ge0} Y_i   1_{\{|Y_i| > 1\}} > x,   \sup_{i \le-1} |Y_i|
\le1  \biggr), \nonumber
\\[-8pt]
\\[-8pt]
\nu^{(u)}(-\infty, -x) &=&  u^{-\alpha}   \Pr \biggl(
u \sum_{i \ge0} Y_i   1_{\{|Y_i| > 1\}} < -x,   \sup_{i \le-1}
|Y_i| \le1  \biggr).
\nonumber
\end{eqnarray}

In case $\alpha\in[1, 2)$, we will need to assume that the
contribution of the smaller increments of the partial sum process is
close to its expectation. The name of the condition is borrowed from
Bartkiewicz et al.~\cite{BaJaMiWi09}, Section~2.4; see Section~\ref{SSSvsv} for a discussion on this assumption.

\begin{cond}[(Vanishing small values)]
\label{cstep6cond} For all $\delta> 0$,
\[
\lim_{u \downarrow0} \limsup_{n \to\infty} \Pr\Biggl [
\max_{0 \le k \le n}  \Biggl| \sum_{i=1}^{k} \biggl ( \frac{X_{i}}{a_{n}}
1_{  \{ \fraca{|X_{i}|}{a_{n}} \le u  \} } - \E \biggl( \frac
{X_{i}}{a_{n}}
1_{  \{ \fraca{|X_{i}|}{a_{n}} \le u  \} }  \biggr)  \biggr)
 \Biggl| > \delta
 \Biggr]=0.
\]
\end{cond}


\begin{theo}
\label{t2} Let $(X_{n})_{n \in\mathbb{N}}$ be a strictly
stationary sequence of random variables, jointly regularly varying
with index $\alpha\in(0,2)$, and of which the tail process
$(Y_{i})_{i \in\ZZ}$ almost surely has no two values of the
opposite sign. Suppose that Conditions~\ref{cfinite-mean-cluster-size} and
\ref{cmixcond} hold. If $1 \le\alpha< 2$, also suppose that
Condition~\ref{cstep6cond} holds. Then the partial sum stochastic
process
\[
V_{n}(t) =
\sum_{k=1}^{[nt]} \frac{X_{k}}{a_{n}} - \lfloor nt \rfloor\E \biggl(
\frac{X_{1}}{a_{n}} 1_{  \{ \fraca{|X_{1}|}{a_{n}} \le1  \} }
 \biggr),  \qquad  t \in[0,1],
\]
satisfies
\[
V_{n} \dto V, \qquad n \to\infty,
\]
in $D[0,1]$ endowed with the $M_{1}$ topology, where $V( \cdot )$
is an $\alpha$-stable L\'{e}vy process with L\'{e}vy triple $(0,
\nu, b)$ given by the limits
\[
\nu^{(u)} \vto\nu, \qquad   \int_{x \dvtx  u < |x| \le1} x   \nu^{(u)}(\rmd
x) - \int_{x \dvtx  u < |x| \le1} x   \mu(\rmd x) \to b
\]
as $u \downarrow0$, with $\nu^{(u)}$ as in \eqref{Enuu} and $\mu$
as in \eqref{Emu}.
\end{theo}

The condition that the tail process cannot switch sign is needed to
ensure continuity of the summation functional; see Lemma~\ref{lprob1}. See Section~\ref{SSSnss} for some discussion of this condition.

\begin{pf}
Note that from Theorem~\ref{Tpointprocesscomplete} and the fact
that $|Y_{n}| \to0$ almost surely as $|n| \to\infty$, the random
variables
\[
u \sum_{j}Z_{ij}1_{\{ |Z_{ij}|>1 \}}
\]
are i.i.d. and almost surely finite. Define
\[
\widehat{N}^{(u)} = \sum_{i} \delta_{(T_{i}^{(u)}, u\sum
_{j}Z_{ij}1_{\{ |Z_{ij}|>1
\}})}.
\]
Then by Proposition 5.3 in Resnick~\cite{Resnick07},
$\widehat{N}^{(u)}$ is a Poisson process (or a~Poisson random
measure) with mean measure
%
\begin{equation}\label{eprodmeas}
\theta u^{-\alpha} \lambda\times F^{(u)},\vadjust{\goodbreak}
\end{equation}
where $\lambda$ is the Lebesgue measure, and $F^{(u)}$ is the
distribution of the random variable $u \sum_{j}Z_{1j}1_{\{
|Z_{1j}|>1 \}}$. But for $0 \le s < t \le1$ and $x>0$, using the
fact that the distribution of $\sum_{j}\delta_{Z_{1j}}$ is equal to
the one of $\sum_{j}\delta_{Y_{j}}$ conditionally on the event $\{
\sup_{i \le-1}|Y_{i}| \le1 \}$, we have
\begin{eqnarray*}
 &&\theta u^{-\alpha} \lambda\times F^{(u)} \bigl([s,t] \times
(x,\infty)\bigr)\\
 && \qquad = \theta u^{-\alpha} (t-s)
F^{(u)}((x,\infty))  \\
&& \qquad  = \theta u^{-\alpha} (t-s)
\Pr\biggl ( u \sum_{j}Z_{1j}1_{\{ |Z_{1j}|>1 \}} >x  \biggr)
\\
&& \qquad  = \theta u^{-\alpha} (t-s) \Pr \biggl( u \sum_{j}Y_{j}1_{\{
|Y_{j}|>1 \}}
>x \Bigm|  \sup_{i \le-1}|Y_{i}| \le1  \biggr)\\
&& \qquad  = \theta u^{-\alpha} (t-s)
\frac{\Pr ( u \sum_{j}Y_{j}1_{\{ |Y_{j}|>1 \}} >x, \sup_{i\le
-1}|Y_{i}| \le1
 )}{\Pr(\sup_{i \le-1}|Y_{i}| \le1)} \\
&& \qquad  = u^{-\alpha} (t-s)
\Pr\biggl ( u \sum_{j}Y_{j}1_{\{ |Y_{j}|>1 \}}
>x, \sup_{i\le-1}|Y_{i}| \le1  \biggr)\\
&& \qquad  = \lambda\times\nu^{(u)}\bigl([s,t] \times(x, \infty)\bigr).
\end{eqnarray*}
The same can be done for the set $[s, t] \times(-\infty, -x)$, so
that the mean measure in \eqref{eprodmeas} is equal to $\lambda
\times\nu^{(u)}$.

Consider now $0<u<v$ and
\[
\psi^{(u)}\bigl (N_{n} | _{[0,1] \times\EE_{u}}\bigr) ( \cdot )
= \psi^{(u)} \bigl(N_{n} | _{[0,1] \times\EE_{v}}\bigr) ( \cdot )
= \sum_{i/n \le  \cdot} \frac{X_{i}}{a_{n}} 1_{  \{ \fraca
{|X_{i}|}{a_{n}} > u
 \} },
\]
which by Lemma~\ref{Lcontsf} converges in distribution in $D[0,1]$
under the $M_{1}$ metric to
\[
\psi^{(u)}\bigl (N^{(v)}\bigr)( \cdot )
=\psi^{(u)} \bigl(N^{(v)} | _{[0,1] \times\EE_{u}}\bigr)( \cdot ).
\]
However, by the definition of the process $N^{(u)}$ in
Theorem~\ref{Tpointprocesscomplete}, it holds that
\[
N^{(u)} \eind
N^{(v)}  |_{[0, 1] \times\EE_u}
\]
for every $v\in(0,u)$.
Therefore the last expression above is equal in distribution to
\[
\psi^{(u)} \bigl(N^{(u)}\bigr)( \cdot )
= \sum_{T_{i}^{(u)} \le  \cdot}
\sum_{j}uZ_{ij}1_{ \{ |Z_{ij}| > 1 \} }.
\]
But since
$\psi^{(u)}(N^{(u)}) = \psi^{(u)} (\widehat{N}^{(u)}) \eind \psi
^{(u)} (\widetilde{N}^{(u)})$,
where
\[
\widetilde{N}^{(u)} = \sum_{i} \delta_{(T_{i}, K_{i}^{(u)})}
\]
is a Poisson process with mean measure $\lambda\times\nu^{(u)}$,
we obtain
\[
\sum_{i = 1}^{\lfloor n   \cdot  \rfloor} \frac{X_{i}}{a_{n}} 1_{
 \{ \fraca{|X_{i}|}{a_{n}} > u
 \} } \dto\sum_{T_{i} \le  \cdot} K_{i}^{(u)}  \qquad\mbox
{as }  n \to\infty
\]
in $D[0,1]$ under the $M_{1}$ metric. From (\ref{eonedimregvar}) we
have, for any $t \in[0,1]$, as $n \to\infty$,
\begin{eqnarray*}
\lfloor nt \rfloor\E \biggl( \frac{X_{1}}{a_{n}}   1_{  \{ u <
\fraca{|X_{1}|}{a_{n}} \le1  \} }
 \biggr) & = & \frac{\lfloor nt \rfloor}{n} \int_{\{x \dvtx  u < |x| \le
1 \}}x
n \Pr \biggl( \frac{X_{1}}{a_{n}} \in\rmd x  \biggr) \\
& \to& t \int_{\{x \dvtx  u < |x| \le1 \}}x   \mu(\rmd x).
\end{eqnarray*}
This convergence is uniform in $t$, and hence
\[
\lfloor n   \cdot  \rfloor\E \biggl( \frac{X_{1}}{a_{n}} 1_{
\{ u < \fraca{|X_{1}|}{a_{n}} \le1  \} }
 \biggr) \to( \cdot ) \int_{\{x \dvtx  u < |x| \le1 \}}x   \mu(\rmd x)
\]
in $D[0,1]$.
Since the latter function is continuous, we can apply
Corollary~12.7.1 in Whitt~\cite{Whitt02}, giving a sufficient
criterion for addition to be continuous. We obtain, as $n \to
\infty$,
%
\begin{eqnarray}
\label{emainconv}
V_{n}^{(u)}( \cdot ) = \sum_{i = 1}^{\lfloor n   \cdot  \rfloor
} \frac{X_{i}}{a_{n}}
1_{  \{ \fraca{|X_{i}|}{a_{n}} > u  \} } - \lfloor n  \cdot
  \rfloor
\E \biggl( \frac{X_{1}}{a_{n}}
1_{  \{ u < \fraca{|X_{1}|}{a_{n}} \le1  \} }
 \biggr) \nonumber
 \\[-8pt]
 \\[-8pt]
\dto V^{(u)}( \cdot ) := \sum_{T_{i} \le  \cdot}
K_{i}^{(u)} - ( \cdot ) \int_{\{x \dvtx  u < |x| \le1 \}}x   \mu
(\rmd x).
\nonumber
\end{eqnarray}
Limit (\ref{emainconv}) can be rewritten as
\begin{eqnarray*}
&&\sum_{T_{i} \le \cdot}
K_{i}^{(u)} - ( \cdot ) \int_{\{x \dvtx  u < |x| \le1 \}}x   \nu
^{(u)}(\rmd x) \\
&& \qquad {}+ ( \cdot ) \biggl ( \int_{\{x \dvtx  u < |x| \le1 \}}x   \nu
^{(u)}(\rmd x)
- \int_{\{x \dvtx  u < |x| \le1 \}}x   \mu(\rmd x)  \biggr).
\end{eqnarray*}
Note that the first two terms represent a L\'{e}vy--Ito
representation of the L\'{e}vy process with characteristic triple
$(0, \nu^{(u)}, 0)$ (see Resnick~\cite{Resnick07}, page~150). The
remaining term is just a linear function of the form $t \mapsto t
b_{u}$. As a consequence, the process $V^{(u)}$ is a L\'{e}vy
process for each $u<1$, with characteristic triple $(0, \nu^{(u)},
b_{u})$, where
\[
b_{u} = \int_{\{x \dvtx  u < |x| \le1 \}}x
\nu^{(u)}(\rmd x) - \int_{\{x \dvtx  u < |x| \le1 \}}x
\mu(\rmd x).
\]

By Theorem 3.1 in Davis and Hsing~\cite{DaHs95}, for $t=1$,
$V^{(u)}(1)$ converges to an $\alpha$-stable random variable. Hence
by Theorem 13.17 in Kallenberg~\cite{Kallenberg97},\vadjust{\goodbreak} there is a
L\'{e}vy process $V( \cdot )$ such that, as $u \to0$,
\[
V^{(u)}( \cdot ) \dto V( \cdot )
\]
in $D[0,1]$ with the $M_{1}$ metric. It has characteristic triple
$(0, \nu, b)$, where $\nu$ is the vague limit of $\nu^{(u)}$ as $u
\to0$ and $b=\lim_{u \to0}b_{u}$ (see Theorem 13.14
in~\cite{Kallenberg97}). Since the random variable $V(1)$ has an
$\alpha$-stable distribution, it follows that the process
$V( \cdot )$ is $\alpha$-stable.

If we show that
\[
\lim_{u \downarrow0} \limsup_{n \to\infty}
\Pr\bigl[d_{M_{1}}\bigl(V_{n}^{(u)}, V_{n}\bigr) > \delta\bigr]=0
\]
for any $\delta>0$, then by Theorem 3.5 in Resnick~\cite{Resnick07} we
will have, as $n \to\infty$,
\[
V_{n} \dto V
\]
in $D[0,1]$ with the $M_{1}$ metric. Since the Skorohod $M_{1}$
metric on $D[0,1]$ is bounded above by the uniform metric on
$D[0,1]$, it suffices to show that
\[
\lim_{u \downarrow0} \limsup_{n \to\infty} \Pr\Bigl (
\sup_{0 \le t \le1} \bigl|V_{n}^{(u)}(t) - V_{n}(t)\bigr| >
\delta \Bigr)=0.
\]
Recalling the definitions, we have
\begin{eqnarray*}
&&\lim_{u \downarrow0}   \limsup_{n \to\infty} \Pr\Bigl (
\sup_{0 \le t \le1} \bigl|V_{n}^{(u)}(t) - V_{n}(t)\bigr| > \delta \Bigr) \\
&& \quad  = \lim_{u \downarrow0} \limsup_{n \to\infty} \Pr \Biggl[
\sup_{0 \le t \le1}  \Biggl| \sum_{i=1}^{\lfloor nt \rfloor} \frac
{X_{i}}{a_{n}}
1_{  \{ \fraca{|X_{i}|}{a_{n}} \le u  \} } - \lfloor nt \rfloor
\E \biggl( \frac{X_{1}}{a_{n}}
1_{  \{ \fraca{|X_{1}|}{a_{n}} \le u  \} }  \biggr) \Biggr | >
\delta
 \Biggr]\\
&& \quad  =\lim_{u \downarrow0} \limsup_{n \to\infty} \Pr \Biggl[
\sup_{0 \le t \le1}  \Biggl| \sum_{i=1}^{\lfloor nt \rfloor}
\biggl\{ \frac{X_{i}}{a_{n}}
1_{  \{ \fraca{|X_{i}|}{a_{n}} \le u  \} } - \E \biggl( \frac
{X_{i}}{a_{n}}
1_{  \{ \fraca{|X_{i}|}{a_{n}} \le u  \} }  \biggr)  \biggr\}
 \Biggr| > \delta
 \Biggr]\\
&& \quad  = \lim_{u \downarrow0} \limsup_{n \to\infty} \Pr \Biggl[
\max_{1 \le k \le n}  \Biggl| \sum_{i=1}^{k}  \biggl\{ \frac{X_{i}}{a_{n}}
1_{  \{ \fraca{|X_{i}|}{a_{n}} \le u  \} } - \E \biggl( \frac
{X_{i}}{a_{n}}
1_{  \{ \fraca{|X_{i}|}{a_{n}} \le u  \} }  \biggr)  \biggr\}
 \Biggr| > \delta
 \Biggr].
\end{eqnarray*}
Therefore we have to show
%
\begin{eqnarray}\label{eslutskycond}
 && \lim_{u \downarrow0} \limsup_{n \to\infty} \Pr \Biggl[
\max_{1 \le k \le n}  \Biggl| \sum_{i=1}^{k}  \biggl\{ \frac{X_{i}}{a_{n}}
1_{  \{ \fraca{|X_{i}|}{a_{n}} \le u  \} } \nonumber\\[-8pt]\\[-8pt]
&&\hphantom{\lim_{u \downarrow0} \limsup_{n \to\infty} \Pr \Biggl[\max_{1 \le k \le n}  \Biggl| \sum_{i=1}^{k}  \biggl\{}- \E \biggl( \frac
{X_{i}}{a_{n}}
1_{  \{ \fraca{|X_{i}|}{a_{n}} \le u  \} }  \biggr)  \biggr\}
 \Biggr| > \delta
 \Biggr]=0.\nonumber
\end{eqnarray}
For $\alpha\in[1,2)$ this relation is simply
Condition~\ref{cstep6cond}. The proof that \eqref{eslutskycond}
automatically holds in case $\alpha\in(0, 1)$ is given at the end of
the proof of Theorem 4.1 in Tyran-Kami{\'n}ska~\cite{TK10SPA}, page 1640.
\end{pf}

\subsection{Discussion of the conditions}
\label{SSdisc}

Here we revisit in detail all the conditions of Theorem~\ref{t2}.

\subsubsection{On joint regular variation}
As we mentioned in the \hyperref[intro]{Introduction}, regular variation of the marginal
distribution with index
$\alpha\in(0,2)$ is both necessary and sufficient for the
existence of an $\alpha$-stable limit for partial sums of i.i.d.
random variables. In Tyran-Kami{\'n}ska~\cite{TK10SPA}, only marginal
regular variation is assumed from the outset, but in combination with
the asymptotic independence condition on the finite-dimensional
distributions, this actually implies joint regular variation.

The joint regular variation assumption \eqref{etailprocess} which
underlies our main result frequently appears in limit theorems for
partial sums~\cite{DaHs95,DaMi98,BaJaMiWi09}. The assumption is
relatively straightforward to verify for many applied models;  see, for
instance,  Section~\ref{Sexamples}. The joint regular variation is
the basis for the point process result of Theorem~\ref{Tpointprocesscomplete}. In particular, it allows us to build on the
theory developed in~\cite{BaSe,DaHs95} to determine the asymptotic
behavior of partial sums over shorter blocks of indices.
On the other hand, we note that 
there are published examples of bounded sequences whose partial sums
have an infinite variance $\alpha$-stable limit (e.g., see Gou\"ezel~\cite{Go04}).

\subsubsection{\texorpdfstring{On the finite mean cluster size Condition \protect\ref{cfinite-mean-cluster-size}}
{On the finite mean cluster size Condition 2.1}}
\label{SSSfmcs}
This assumption, which appears frequently in the literature
\cite{BaSe,DaHs95,Segers03,Segers05,Smith92}, restricts the length of
clusters of extremes. It implies that the (max)-stable attractors of
appropriately normalized partial sums and maxima have the same index
$\alpha$ as the ones for the associated i.i.d. sequence. Alternative
assumptions of this kind also exist, most of which are stronger; see,
for instance,~\cite{BaJaMiWi09} for a~short review.


\subsubsection{\texorpdfstring{On the $\mathcal{A}(a_n)$ mixing Condition~\protect\ref{cmixcond}}
{On the A(a n) mixing Condition 2.2}}\label{SSSmixing}

Extremely rich literature exists on mixing conditions and their
relation with limit theorems. Our assumption $\mathcal{A}'(a_{n})$ is
a recognizable extension of the mixing condition $\mathcal{A}(a_{n})$
due to~\cite{DaHs95}. Like the latter condition it is implied by the
more frequently used strong mixing property (see~\cite{Kr10}).
However, if one is only interested in the limiting behavior of partial
sums, weaker assumptions suffice (see~\cite{BaJaMiWi09}).

\subsubsection{\texorpdfstring{On the vanishing small values Condition~\protect\ref{cstep6cond}}
{On the vanishing small values Condition 3.3}}\label{SSSvsv}

The name of the condition is borrowed from~\cite{BaJaMiWi09}. Similar
conditions are ubiquitous in the related literature on the limit theory
for partial sums~\cite{Avram92,Durrett78,Leadbetter88,TK10SPA}. In
case $\alpha\in(0,1)$, it is simply a consequence of regular
variation, and in the i.i.d. case, it also holds for $\alpha\in[1,
2)$ (see Resnick~\cite{Resnick86}). More generally, Tyran-Kami{\'n}ska~\cite{TK10SPA} showed that the condition holds if the sequence
has $\rho$-mixing coefficients which satisfy $\sum_{j\geq1} \rho
(2^j)<\infty$. For linear processes of which the coefficients decay
sufficiently fast, Tyran-Kami{\'n}ska~\cite{TK10SPL} showed that the
condition can be omitted.

\subsubsection{{About the no sign switching condition}}\label{SSSnss}

The assumption that the tail process has no two values of the opposite
sign is crucial to obtain weak convergence of the partial sum
process
in the $M_{1}$ topology. It is admittedly restrictive but unavoidable
since the $M_1$ topology, roughly speaking, can handle several
(asymptotically) instantaneous jumps only if they are in the same
direction (see Avram and Taqqu~\cite{Avram92}, Section 1, and
Whitt~\cite{Whitt02}, Chapter~12). Note that unlike Dabrowski and
Jakubowski~\cite{DaJa}, our assumption does not exclude nonassociated
sequences in general because it involves only the tail dependence in
the process.

In Avram and Taqqu~\cite{Avram92}, Section 1, a conjecture is
formulated concerning convergence in Skorohod's $M_2$ topology, which
is somewhat weaker than the $M_1$ topology. Rather than being all of
the same sign, extremes values within a cluster should be such that the
values of the partial sums during the cluster are all contained in the
interval formed by the partial sums at the beginning and the end of a cluster.


 There appear to be some ways of omitting the no sign
switching condition altogether. Neither of them is pursued here,
however. First, one could opt for a much weaker topology on
$D[0,1]$, like $L_{1}$, for instance. Another possibility is to
avoid the within-cluster fluctuations in the partial sum process,
for example,
  by smoothing out its trajectories or by considering the
process $t \mapsto S_{r_{n} \lfloor k_{n}t \rfloor}$. If we do so,
then convergence actually holds in the stronger $J_1$ topology
(see Krizmani\'c ~\cite{Kr10}, Chapter 3).

%

\subsection{Simplifications}
\label{SSsimpl}

In certain cases, the formula for the L\'evy measure can be simplified.
Moreover, if $\alpha\in(0, 1)$, then no centering is needed.

\subsubsection{A closed form expression for the limiting L\'evy measure}\label{RLevyChar3}

It turns out that if the spectral tail process $(\Theta_i)_{i \in
\ZZ}$ satisfies an additional integrability condition, the formula
for the L\'evy measure $\nu$ simplifies considerably. Note that
in our case the L\'evy measure $\nu$ satisfies the scaling
property
\[
\nu(s   \cdot ) = s^{-\alpha} \nu(  \cdot )
\]
(see Theorem 14.3 in Sato~\cite{Sato99}). In particular, $\nu$ can be
written as
\[
\nu(\rmd x)
=  \bigl( c_+   1_{(0, \infty)}(x) + c_-   1_{(-\infty,0)}(x)
\bigr)   \alpha|x|^{-\alpha-1} \, \rmd x
\]
for some nonnegative constants $c_+$ and $c_-$, and therefore
$\nu(\{x \})=0$ for every $x \in\EE$. Thus, from Theorem 3.2 in
Resnick~\cite{Resnick07} we have
\begin{eqnarray*}
c_+
&=& \nu(1,\infty) \\
&=& \lim_{u \to0} \nu^{(u)}(1,\infty) \\
&=& \lim_{u \to0} u^{-\alpha}   \Pr \biggl( u \sum_{i \ge0} Y_i
  1_{\{|Y_i| > 1\}} > 1,   \sup_{i \le-1} |Y_i| \le1  \biggr) \\
&=& \lim_{u \to0} u^{-\alpha}   \int_1^\infty\Pr \biggl( u \sum
_{i \ge0} r \Theta_i   1_{\{r|\Theta_i| > 1\}} > 1,   \sup_{i \le
-1} r|\Theta_i| \le1  \biggr)  \,\rmd(-r^{-\alpha}) \\[-2pt]
&=& \lim_{u \to0} \int_{u}^{\infty}
\Pr \biggl( \sum_{i \ge0} r \Theta_{j}   1_{\{ r|\Theta_{j}|>u \}
} > 1,   \sup_{i \le-1} r|\Theta_{i}| \le u  \biggr)\,\rmd
(-r^{-\alpha}),
\end{eqnarray*}
and similarly
\[
c_-
= \lim_{u \to0} \int_{u}^{\infty}
\Pr \biggl( \sum_{i \ge0} r \Theta_{j}   1_{\{ r|\Theta_{j}|>u \}}
< -1, \sup_{i \le-1}r|\Theta_{i}| \le u  \biggr)
\, \rmd(-r^{-\alpha}).
\]

Now suppose further that
%
\begin{equation}
\label{Esumfinite}
\E \biggl[ \biggl ( \som_{i \ge0} |\Theta_i|  \biggr)^\alpha \biggr] <
\infty.
\end{equation}
Then by the dominated convergence theorem,
%
\begin{eqnarray}
\label{Ecplus}
c_+
&=& \int_0^\infty\Pr\biggl ( \sum_{i \ge0} r \Theta_i > 1 ;
\forall i \le-1 \dvtx  \Theta_i = 0  \biggr) \, \rmd(-r^{-\alpha}) \nonumber
\\[-9pt]
\\[-9pt]
&=& \E \biggl[ \biggl \{ \max \biggl( \som_{i \ge0} \Theta_i, 0
\biggr)  \biggr\}^\alpha  1_{\{\forall i \le-1 \dvtx  \Theta_i = 0\}}  \biggr],
\nonumber\\[-2pt]
\label{Ecminus}
c_-
&=& \E\biggl [  \biggl\{ \max \biggl( - \som_{i \ge0} \Theta_i, 0
 \biggr)  \biggr\}^\alpha  1_{\{\forall i \le-1 \dvtx  \Theta_i = 0\}}
 \biggr].
\end{eqnarray}
These relations can be applied to obtain the L\'evy measure $\nu$
for certain heavy-tailed moving average processes
(Example~\ref{exFiniteMA}).

\subsubsection{About centering}\label{rcent}

If $\alpha\in(0,1)$, the centering function in
the definition of the stochastic process $V_{n}( \cdot )$ can be
removed. This affects the characteristic triple of the limiting
process in the way we describe here.

By Karamata's theorem, as $n \to
\infty$,
\[
n \E \biggl( \frac{X_{1}}{a_{n}}   1_{  \{ \fraca{|X_{1}|}{a_{n}}
\le1  \} }
 \biggr) \to(p-q) \frac{\alpha}{1-\alpha}
\]
with $p$ and $q$ as in \eqref{Emu}. Thus, as $n \to\infty$,
\[
\lfloor n   \cdot  \rfloor\E \biggl( \frac{X_{1}}{a_{n}}   1_{
 \{ \fraca{|X_{1}|}{a_{n}} \le1  \} }
 \biggr) \to( \cdot )(p-q) \frac{\alpha}{1-\alpha}
\]
in $D[0,1]$, which leads to
\[
\sum_{k=1}^{\lfloor n   \cdot  \rfloor} \frac{X_{k}}{a_{n}}
\dto V( \cdot ) + ( \cdot ) (p-q)\frac{\alpha}{1-\alpha}
\]
in $D[0,1]$ endowed with the $M_{1}$ topology. The characteristic
triple of the limiting process is therefore $(0, \nu, b')$ with $b' =
b + (p-q)\alpha/ (1-\alpha)$.\vadjust{\goodbreak}

\section{Examples}\label{Sexamples}

In case of asymptotic independence, the limiting stable L\'evy
process is the same as in the case of an i.i.d. sequence with the
same marginal distribution (Examples~\ref{exDcond} and
\ref{exStochVol}). Heavy-tailed moving averages and $\operatorname{GARCH}(1,1)$
processes (Examples~\ref{exFiniteMA} and~\ref{exGARCH},
respectively) yield more interesting limits.

\begin{ex}[(Isolated extremes models)]
\label{exDcond} Suppose $(X_{n})$ is a strictly stationary and
strongly mixing sequence of regularly varying random variables with
index $\alpha\in(0,2)$ that satisfies the dependence condition
$D'$ in Davis~\cite{Da83}, that is,
\[
\lim_{k \to\infty} \limsup_{n \to\infty} n \sum_{i=1}^{\lfloor n/k
\rfloor} \Pr \biggl( \frac{|X_{0}|}{a_{n}} > x,
\frac{|X_{i}|}{a_{n}} > x  \biggr) = 0 \qquad\mbox{for all }  x
>0,
\]
where $(a_{n})_n$ is a positive real sequence such that $n
\Pr(|X_{0}|>a_{n}) \to1$ as $n \to\infty$. Condition $D'$ implies
\[
\Pr(|X_i| > a_n \mid|X_0| > a_n) = \frac{n \Pr(|X_0| > a_n,
|X_i| > a_n)}{n \Pr(|X_0| > a_n)} \to0  \qquad\mbox{as $n \to
\infty$}
\]
for all positive integer $i$; that is, the variables $|X_0|$ and
$|X_i|$ are asymptotically independent. As a consequence, the series
$(X_n)_n$ is regularly varying and its tail process is the same as
that for an i.i.d. sequence; that is, $Y_n = 0$ for $n \ne0$, and
$Y_0$ is as described in Section~\ref{SSstatpointtail}. It is
trivially satisfied that no two values of $(Y_{n})_n$ are of the
opposite sign.

Since the sequence $(X_{n})$ is strongly mixing, Condition~\ref{cmixcond} is verified.
The finite mean cluster size Condition~\ref{cfinite-mean-cluster-size}
follows from condition $D'$, for the latter implies
\[
\lim_{n \to\infty} n \sum_{i=1}^{r_{n}} \Pr \biggl( \frac
{|X_{0}|}{a_{n}} > x,
\frac{|X_{i}|}{a_{n}} > x  \biggr) = 0 \qquad\mbox{for all }  x
>0
\]
for any positive integer sequence $(r_{n})_n$ such that $r_{n} \to
\infty$ and $r_{n} / n \to0$ as \mbox{$n \to\infty$}.

If we additionally assume that the sequence $(X_{n})$ satisfies the
vanishing small values Condition~\ref{cstep6cond} in case $\alpha
\in[1,2)$, then by Theorem~\ref{t2} the sequence of partial sum
stochastic processes $V_{n}( \cdot )$ converges in $D[0,1]$ with the
$M_{1}$ topology to an $\alpha$-stable L\'{e}vy process $V( \cdot
)$ with characteristic triple $(0, \mu, 0)$ with $\mu$ as in \eqref{Emu}, just as in the i.i.d. case. It can be shown that the above
convergence holds also in the $J_{1}$ topology (see Krizmani\'{c}~\cite{Kr10}).

Condition~\ref{cstep6cond} applies, for instance, if the series
$(X_{n})_n$ is a function of a~Gaussian causal ARMA process, that is,
$X_{n} = f(A_{n})$, for some Borel function $f \dvtx  \mathbb{R} \to
\mathbb{R}$ and some Gaussian causal ARMA process $(A_{n})_n$. From
the results in Brockwell and Davis~\cite{BrDa91} and Pham and
Tran~\cite{PhTr85} (see also Davis and Mikosch~\cite{DaMi09}) it\vadjust{\goodbreak}
follows that the sequence $(A_{n})_n$ satisfies the strong mixing
condition with geometric rate. In this particular case this implies
that the sequence $(A_{n})_n$ satisfies the $\rho$-mixing condition
with geometric rate (see Kolmogorov and Rozanov
\cite{KoRo60}, Theorem 2), a property which transfers immediately to
the series $(X_{n})_n$. Hence by Tyran-Kami{\'n}ska~\cite{TK10SPA}, Lemma 4.8, the vanishing small values Condition~\ref{cstep6cond} holds.
\end{ex}

\begin{ex}[(Stochastic volatility models)]
\label{exStochVol} Consider the stochastic volatility model
\[
X_{n} = \sigma_{n} Z_{n}, \qquad n \in\mathbb{Z},
\]
where the noise sequence $(Z_{n})$ consists of i.i.d. regularly
varying random variables with index $\alpha\in(0,2)$, whereas the
volatility sequence $(\sigma_{n})_n$ is strictly stationary, is
independent of the sequence $(Z_{n})_n$ and consists of positive
random variables with finite moment of the order $4 \alpha$.

Since the random variables $Z_{i}$ are independent and regularly
varying, it follows that the sequence $(Z_{n})_n$ is regularly
varying with index $\alpha$. By an application of the multivariate
version of Breiman's lemma, the sequence $(X_{n})_n$ is regularly
varying with
index $\alpha$ too{; see Basrak et al.~\cite{BDM02b}, Proposition~5.1}.

From the results in Davis and Mikosch~\cite{DaMi08}, it follows that
%
\begin{equation}
\label{eanticl1}
n \sum_{i=1}^{r_{n}} \Pr(|X_{i}|>ta_{n}, |X_{0}|>ta_{n}) \to0
\qquad\mbox{as $n \to\infty$}
\end{equation}
for any $t>0$, where $(r_{n})_n$ is a sequence of positive integers
such that $r_{n} \to\infty$ and $r_{n} / n \to0$ as $n \to\infty$,
and $(a_{n})_n$ is a positive real sequence such that $n \Pr
(|X_{1}|>a_{n}) \to1$ as $n \to\infty$. From this relation, as in
Example~\ref{exDcond}, it follows that the finite mean cluster size
Condition~\ref{cfinite-mean-cluster-size} holds. Moreover, the tail
process $(Y_n)_n$ is the same as in the case of an i.i.d. sequence,
that is, $Y_n = 0$ for $n \ne0$. In particular, the tail process has
no two values of the opposite sign.

Assume that $(\log\sigma_{n})_n$ is a Gaussian casual ARMA process.
Then $(X_{n})_n$ satisfies the strong mixing condition with geometric
rate (see Davis and Mikosch~\cite{DaMi09}). Hence the $\mathcal
{A}'(a_n)$ mixing Condition~\ref{cmixcond} holds.

In case $\alpha\in[1, 2)$, we also assume the vanishing small values
Condition~\ref{cstep6cond} holds. Then all conditions in
Theorem~\ref{t2} are satisfied, and we obtain the convergence of the
partial sum stochastic process toward an $\alpha$-stable L\'{e}vy
process with characteristic triple $(0,\mu,0)$, with $\mu$ as
in~\eqref{Emu}.
\end{ex}

\begin{ex}[({Moving averages})]
\label{exFiniteMA}
Consider the finite-order moving average defined
by
\[
X_{n} = \sum_{i=0}^{m} c_{i}Z_{n-i}, \qquad n \in\mathbb{Z},
\]
where $(Z_{i})_{i\in\mathbb{Z}}$ is an i.i.d. sequence of
regularly varying random variables with index $\alpha\in(0,2)$, $m
\in\mathbb{N}$, $c_{0}, \ldots, c_{m}$ are\vadjust{\goodbreak} nonnegative constants
and at least $c_{0}$ and $c_{m}$ are not equal to $0$. Take a
sequence of positive real numbers $(a_{n})$ such that
%
\begin{equation}\label{enizan}
n \Pr(|Z_{1}|>a_{n}) \to1 \qquad\mbox{as }  n \to\infty.
\end{equation}

The finite-dimensional distributions of the series $(X_n)_n$ can be
seen to be multivariate regularly varying by an application of
Proposition 5.1 in Basrak et al.~\cite{BDM02b} (see also Davis and
Resnick~\cite{DaRe85}). Moreover, if we assume (without loss of
generality) that $\sum_{i=0}^{m} c_{i}^{\alpha} = 1$, then also
\[
n \Pr(|X_{0}|>a_{n}) \to1 \qquad\mbox{as }  n \to\infty.
\]
The tail process $(Y_n)_n$ in \eqref{etailprocess} of the series
$(X_n)_n$ can be found by direct calculation (see also Meinguet and
Segers~\cite{MS10}, Proposition 8.1, for an extension to infinite-order
moving averages). First, $Y_0 = |Y_0|
\Theta_0$ where $|Y_0|$ and $\Theta_0 = \operatorname{sign}(Y_0)$
are independent with $P(|Y_0| > y) = y^{-\alpha}$ for $y \ge1$ and
$\Pr(\Theta_0 = 1) = p = 1 - \Pr(\Theta_0 = -1)$. Next, let $K$
denote a random variable with values in the set $\{0,\ldots, m\}$,
independent of $Y_0$ and such that $\Pr(K=k)= |c_k|^{\alpha}$
(recall the assumption $\sum_{i=0}^{m} c_{i}^{\alpha}=1$). To
simplify notation, put $c_i :=0 $ for $i \notin\{ 0,\ldots, m \}$.
Then
\[
Y_n = (c_{n+K}/c_K)   Y_0, \qquad\Theta_n = (c_{n+K}/c_K)   \Theta
_0, \qquad n \in\ZZ,
\]
represents the tail process and spectral process of $(X_n)_n$,
respectively. Clearly, at most $m+1$ values $Y_n$ and $\Theta_n$ are
different from 0 and all have the same sign.

Since the sequence $(X_{n})_n$ is $m$-dependent, it is also strongly
mixing, and therefore the $\mathcal{A}'(a_n)$ mixing Condition~\ref{cmixcond} holds. By the same property it is easy to see that the
finite mean cluster size Condition~\ref{cfinite-mean-cluster-size}
holds. Moreover, in view of Lemma~4.8 in Tyran-Kami{\'n}ska~\cite{TK10SPA}, the vanishing small values Condition~\ref{cstep6cond}
holds as well when $\alpha\in[1,2)$.

As a consequence, the sequence {$(X_n)_n$} satisfies all the
conditions of Theorem~\ref{t2}, and the partial sum process
converges toward a stable L\'evy process~$V( \cdot )$. The
L\'evy measure $\nu$ can be derived from Section~\ref{RLevyChar3}:
since \eqref{Esumfinite} is trivially fulfilled, we obtain from
\eqref{Ecplus} and \eqref{Ecminus},
\[
\nu(\rmd x) =  \Biggl( \som_{i=0}^{m} c_{i}  \Biggr)^{\alpha}
\bigl( p   1_{(0,\infty)}(x) + q   1_{(-\infty,0)}(x)  \bigr)   \alpha
|x|^{-1-\alpha} \, \rmd x,
\]
which corresponds with the results in Davis and
Resnick~\cite{DaRe85} and Davis and Hsing~\cite{DaHs95}. Further, if
$\alpha\in(0,1) \cup(1,2)$, then in the latter two references it
is shown that
\[
b = (p-q)   \frac{\alpha}{1-\alpha}   \Biggl \{  \Biggl( \som
_{i=0}^{m} c_{i}  \Biggr)^{\alpha} -1  \Biggr\},
\]
with $q = 1-p$. The case when $\alpha=1$ can be treated similarly, but
the corresponding expressions are much more complicated and are omitted
here (see Davis and Hsing~\cite{DaHs95}, Theorem 3.2 and Remark 3.3).\vadjust{\goodbreak}

Infinite-order moving averages with nonnegative coefficients are
considered in Avram and Taqqu~\cite{Avram92} and Tyran-Kami{\'n}ska~\cite{TK10SPL}.
The idea is to approximate such processes by a
sequence of finite-order moving averages, for which Theorem~\ref{t2}
applies, and to show that the error of approximation is negligible in
the limit.
\end{ex}

\begin{ex}[(ARCH/GARCH models)]
\label{exGARCH} We consider the $\operatorname{GARCH}(1,1)$ model
\[
X_{n}=\sigma_{n} Z_{n},
\]
where $(Z_{n})_{n \in\mathbb{Z}}$ is a sequence of i.i.d. random
variables with $\E(Z_{1}) = 0$ and $\operatorname{var}(Z_{1}) = 1$,
and
%
\begin{equation}\label{estochvol}
\sigma_{n}^{2} = \alpha_{0} + (\alpha_{1} Z_{n-1}^{2} +
\beta_{1}) \sigma_{n-1}^{2},
\end{equation}
with $\alpha_{0}, \alpha_{1}, \beta_{1}$ being nonnegative
constants. Assume that $\alpha_{0}>0$ and
\[
-\infty\le\E\ln(\alpha_{1}Z_{1}^{2} + \beta_{1}) < 0.
\]
Then there exists a strictly stationary solution to the stochastic
recurrence equation \eqref{estochvol} (see Goldie~\cite{Goldie91}
and Mikosch and St\u{a}ric\u{a}~\cite{MiSt00}). The process $(X_{n})$
is then strictly stationary too. If $\alpha_{1}>0$ and $\beta_{1}>0$
it is called a~$\operatorname{GARCH}(1,1)$ process, while if $\alpha_{1}>0$ and
$\beta_{1}=0$ it is called an ARCH(1) process.

In the rest of the example we consider a stationary squared
$\operatorname{GARCH}(1,1)$ process $(X_{n}^{2})_n$.
Assume that $Z_{1}$ is symmetric, has a positive
Lebesgue density on $\mathbb{R}$ and there exists $\alpha\in(0,2)$ such
that
\[
\E[(\alpha_{1}Z_{1}^{2} + \beta_{1})^{\alpha}]=1
\quad\mbox{and} \quad
\E[(\alpha_{1}Z_{1}^{2} + \beta_{1})^{\alpha} \ln(\alpha
_{1}Z_{1}^{2} + \beta_{1})] < \infty.
\]
Then it is known that the processes $(\sigma_{n}^{2})_n$ and
$(X_{n}^{2})_n$ are regularly varying with index $\alpha$ and
strongly mixing with geometric rate~\cite{BDM02b,MiSt00}. Therefore
the sequence $(X_{n}^{2})_n$ satisfies the $\mathcal{A}'(a_n)$ mixing
Condition~\ref{cmixcond}. The finite mean cluster size Condition~\ref{cfinite-mean-cluster-size} for the sequence $(X_{n}^{2})_n$
follows immediately from the results in Basrak et al.~\cite{BDM02b}.

The (forward) tail process of the bivariate sequence
$((\sigma_{n}^{2},X_{n}^{2}))_n$ is not too difficult to
characterize (see Basrak and Segers~\cite{BaSe}). Obviously, the tail
process of $(X_{n}^{2})_n$ cannot have two values of the opposite
sign.

If additionally the vanishing small values Condition~\ref{cstep6cond} holds when $\alpha\in
[1,2)$, then by Theorem~\ref{t2}, the sequence of partial sum
stochastic processes $(V_{n}( \cdot ))_n$, defined by
\[
V_{n}(t)
= \sum_{k=1}^{[nt]} \frac{X_{k}^{2}}{a_{n}} - \lfloor nt \rfloor \E
 \biggl(
\frac{X_{1}^{2}}{a_{n}} 1_{  \{ \fraca{X_{1}^{2}}{a_{n}} \le1
 \} }  \biggr), \qquad t \in[0,1],
\]
converges weakly to an $\alpha$-stable L\'{e}vy process
$V( \cdot )$ in $D[0,1]$ under the~$M_{1}$ topology. Here
$(a_n)_n$ is a positive sequence such that $n \Pr(X_0^2 > a_n) \to
1$ as $n \to\infty$.\vadjust{\goodbreak}

In case $\alpha\in(0,1) \cup(1,2)$, the characteristic triple
$(0, \nu, b)$ of the stable random variable $V(1)$ and thus of the
stable L\'evy process $V( \cdot )$ can be determined from
Bartkiewicz et al.~\cite{BaJaMiWi09}, Proposition 4.8, Davis and
Hsing~\cite{DaHs95}, Remark 3.1, and Section~\ref{rcent}: after
some calculations, we find
\[
\nu(\rmd x) = c_+   1_{(0, \infty)}(x)   \alpha x^{-\alpha-1} \,\rmd x,
 \qquad
b = \frac{\alpha}{1-\alpha} (c_{+}-1),
\]
where
\[
c_{+}  = \frac{\E[(Z_{0}^{2} + \widetilde{T}_{\infty})^{\alpha} -
\widetilde{T}_{\infty}^{\alpha}]}{ \E(|Z_{1}|^{2\alpha})},  \qquad
\widetilde{T}_{\infty}  = \sum_{t=1}^{\infty} Z_{t+1}^{2} \prod
_{i=1}^{t} (\alpha_{1} Z_{i}^{2} + \beta_{1}).
\]
\end{ex}

\section*{Acknowledgments}
The authors wish to acknowledge helpful comments by the
Editor and the reviewer, in particular those in connection with
the literature survey. Stimulating discussions with participants
of the Conference on Latest Developments in Heavy-Tailed
Distributions (Universit\'e libre de Bruxelles, March 26--27,
2010) and the seminar series at SAMOS (Universit\'e Paris I, June
18, 2010) are greatly appreciated.


%

\printaddresses

\end{document}